\documentclass{amsart}

\usepackage{xypic}
\xyoption{all}
\usepackage{epsfig}
\usepackage{amsthm}
\usepackage{amssymb}
\usepackage{amsmath}
\usepackage{amscd}

\newtheorem {Theorem}  {Theorem}

\numberwithin{Theorem}{section}

\newtheorem {Lemma}[Theorem]  {Lemma}
\newtheorem {Proposition}[Theorem]{Proposition}
\theoremstyle{definition}

\theoremstyle{remark}

\newtheorem {Corollary}[Theorem]{Corollary}

\expandafter\chardef\csname pre amssym.def
at\endcsname=\the\catcode`\@ \catcode`\@=11
\def\undefine#1{\let#1\undefined}
\def\newsymbol#1#2#3#4#5{\let\next@\relax
 \ifnum#2=\@ne\let\next@\msafam@\else
 \ifnum#2=\tw@\let\next@\msbfam@\fi\fi
 \mathchardef#1="#3\next@#4#5}
\def\mathhexbox@#1#2#3{\relax
 \ifmmode\mathpalette{}{\m@th\mathchar"#1#2#3}%
 \else\leavevmode\hbox{$\m@th\mathchar"#1#2#3$}\fi}
\def\hexnumber@#1{\ifcase#1 0\or 1\or 2\or 3\or 4\or 5\or 6\or 7\or 8\or
 9\or A\or B\or C\or D\or E\or F\fi}

\font\teneufm=eufm10 \font\seveneufm=eufm7 \font\fiveeufm=eufm5
\newfam\eufmfam
\textfont\eufmfam=\teneufm \scriptfont\eufmfam=\seveneufm
\scriptscriptfont\eufmfam=\fiveeufm

\catcode`\@=\csname pre amssym.def at\endcsname

\newcounter{remark}
\newenvironment{remark}
{\medskip \stepcounter{remark} \noindent \textit{Remark
\arabic{section}.\arabic{remark}.}}{\rm \cbdu}

\def\3n{\negthinspace \negthinspace \negthinspace }
\def\2n{\negthinspace \negthinspace }
\def\1n{\negthinspace }

\newcommand{\bb}{\bibitem}
\newcommand{\bg}{\begin{equation}}
\newcommand{\ed}{\end{equation}}
\newcommand{\bga}{\begin{eqnarray}}
\newcommand{\eda}{\end{eqnarray}}
\newcommand{\pf}{\textbf{Proof \ }}

\def\cbdu{\hfill{$\Box$}}

\renewcommand{\a}{\alpha}
\renewcommand{\th}{\theta}

\newcommand{\R}{\mathbf{R}}

\newcommand{\Ff}{{\mathcal F}}

\newcommand{\eledos}{L^2}
\newcommand{\eledosdelta}{L^{\frac{m}{\delta}}}
\newcommand{\eleq}{L^q}
\newcommand{\nonlinear}{u \cdot \nabla \th}
\newcommand{\fourierzero}{\widehat{\theta_0}}
\newcommand{\kernel}{K_{\a}}
\newcommand{\nonlinearapprox}{(u_n \cdot \nabla) \th_n}
\newcommand{\hattheta}{\hat{\theta}}
\newcommand{\hatpsi}{\widehat{\psi}}
\newcommand{\exponent}{\frac{2}{2\a -1}}

\def  \R   {{\mathbb R}}

\def  \12  {{\frac{1}{2}}}

\def\bd{\begin{definition}}
\def\ede{\end{definition}}
\def\be{\begin{equation}}
\def\bel{\begin{equation}\label}
\def\ee{\end{equation}}
\def\bt{\begin{Theorem}}
\def\et{\end{Theorem}}
\def\bc{\begin{Corollary}}
\def\ec{\end{Corollary}}
\def\bl{\begin{Lemma}}
\def\el{\end{Lemma}}
\def\bp{\begin{Proposition}}
\def\ep{\end{Proposition}}

\def\br{\begin{remark}}
\def\er{\end{remark}}
\def\ba{\begin{array}}
\def\ea{\end{array}}

\def\bea{\begin{eqnarray}}
\def\eea{\end{eqnarray}}

\begin{document}






\title[Decay of solutions to quasi-geostrophic equation]{Decay of weak solutions to the 2D dissipative quasi-geostrophic equation}

\author[C.J. Niche and M.E. Schonbek]{C\'esar J. Niche and Mar\'{\i}a E. Schonbek}

\address{Department of Mathematics, UC Santa Cruz, Santa Cruz, CA 95064, USA}

\email{cniche@math.ucsc.edu}
\email{schonbek@math.ucsc.edu}

\thanks{The second author was partially supported by NSF grant DMS-0600692}

\begin{abstract}
We address the decay of the norm of weak solutions to the 2D dissipative quasi-geostrophic equation. When the initial data $\theta_0$ is in $L^2$ only, we prove that the $L^2$ norm tends to zero but with no uniform rate, that is, there are solutions with arbitrarily slow decay. For $\theta_0$ in $L^p \cap \eledos$, with $1 \leq p < 2$, we are able to obtain a uniform decay rate in $L^2$. We also prove that when the $L^{\exponent}$ norm of $\theta_0$ is small enough, the $L^q$ norms, for $q > \exponent$, have uniform decay rates. This result allows us to prove decay for the $L^q$ norms, for $q \geq \exponent$, when $\theta_0$ is in $\eledos \cap L^{\exponent}$. 
\end{abstract}

\maketitle

\section{Introduction and statement of results}

We consider the dissipative 2D quasi-geostrophic equation

\bga \label{eq:qge-two}
\theta_{t} + (u \cdot \nabla) \theta + (- \Delta)^{\a} \theta = 0   \nonumber \\ \th (x,0) = \th_0(x)
\eda
where $x \in \R ^2, t > 0$ and super-critical exponent $\frac{1}{2} < \a \leq 1 $. In this equation, $\theta = \theta(x,t)$ is a real scalar function {\em (the temperature of the fluid)}, $u$ is an incompressible vector field {\em (the velocity of the fluid)} determined by the scalar function $\psi$ ({\em the stream function)} through

\begin{displaymath}
u = (u_1, u_2) = (-\frac{\partial \psi}{\partial x_2}, \frac{\partial \psi}{\partial x_1}).
\end{displaymath}
The temperature $\theta$ and the stream function $\psi$ are related by

\begin{displaymath}
\Lambda \psi = - \theta
\end{displaymath}
where $\Lambda$ is the usual operator given by $\Lambda = (- \Delta)^{\frac{1}{2}}$ and defined via the Fourier transform as

\begin{displaymath}
\widehat{\Lambda^{\frac{s}{2}} f} (\xi) = |\xi| ^s \hat{f} (\xi), \qquad s \geq 0.
\end{displaymath}

When $\a = \frac{1}{2}$, ``dimensionally, the 2D quasi-geostrophic equation is the analogue of the 3D Navier-Stokes equations'' (Constantin and Wu \cite{constantin-wu}), and the behaviour of solutions to (\ref{eq:qge-two}) is similar to that of the 3D Navier-Stokes equations. For this reason, $\a = \frac{1}{2}$ is considered the {\em critical exponent}, while $\a \in (\frac{1}{2},1]$ are the {\em supercritical exponents}. Note that when $\a = 1$, (\ref{eq:qge-two}) is the vorticity equation of the 2D Navier-Stokes equations. Besides its intrinsic mathematical interest, the dissipative 2D quasi-geostrophic equation describes models arising in meteorology and oceanography. More specifically, it can be derived from the General Quasi Geostrophic equations by assuming constant potential vorticity and constant buoyancy frequency (see Constantin, Majda and Tabak \cite{const-majda-tabak} and Pedlosky \cite{pedlosky}).

Consider the dissipative quasi-geostrophic equation with supercritical exponent, this is $\a \in (\frac{1}{2},1]$. In this article, we address the uniform decay of the $L^q$ norm, for $q \geq 2$,  of weak solutions to (\ref{eq:qge-two}) for the initial data $\theta _0$ in different spaces. We first describe results related to the ones obtained here. 

In his Ph.D. thesis, Resnick \cite{resnick} proved existence of global solutions to (\ref{eq:qge-two}) for $\theta_0$ in $L^2$. Moreover, he proved a maximum principle 

\be \label{eq:max-principle}
\| \theta(t) \| _{L^p} \leq \| \theta_0 \| _{L^p}, \quad t \geq 0
\ee
for $1 < p \leq \infty$. Constantin and Wu \cite{constantin-wu} established uniqueness of ``strong'' solutions (for a precise statement of this and Resnick's result, see Section \ref{basic-facts}) and also showed that for $\theta_0$ in $L^2 \cap L^1$

\be \label{eq:decay-constantin-wu} 
\Vert \theta (t) \Vert _{L^2} \leq C (1 + t)^{-\frac{1}{2 \a}}, \quad t \geq 0.
\ee
Their proof relies on an adaptation of the Fourier splitting method developed by Schonbek \cite{sch-2}, \cite{sch-3} and on the retarded mollifiers method of Cafarelli, Kohn and Nirenberg \cite{ckn}. Moreover, they proved that for generic initial data, the decay rate (\ref{eq:decay-constantin-wu}) is optimal. Using rather general pointwise estimates for the fractional derivative $\Lambda ^{\a} \theta$ and a positivity lemma, C\'ordoba and C\'ordoba \cite{cor-cor} gave a new proof of (\ref{eq:max-principle}) and proved decay of solutions when $\theta _0$ is in $L^1 \cap L^p$, for $ 1 < p < \infty$. More specifically, they showed that

\be \label{eq:bound-cordoba-cordoba} 
\| \theta(t) \| _{L^p} \leq C_1 \left( 1 + C_2 t \right) ^{-\frac{p-1}{\a p}}, \quad t \geq 0
\ee
where $C_1$ and $C_2$ are explicit constants. Working along the same lines, Ju \cite{ju-2} obtained an improved maximum principle of the form

\be \label{eq:max-prin-ju} 
\| \theta(t) \| _{L^p} \leq \| \theta_0  \| _{L^p} \left( 1 + \frac{C^{\frac{1}{p-2}}}{p - 2} \, t \right) ^{\frac{2-p}{2 p \a}} 
\ee 
for $\theta_0$ in $L^2 \cap L^p$, with $p \geq 2$ and a constant $C \neq 1$. Note that for $p = 2$, i.e. $\theta_0$ in $L^2$, this expression reduces to (\ref{eq:max-principle}), this is  $\| \theta(t) \| _{\eledos} \leq \| \theta_0 \| _{\eledos}$. \\

We now state the results we prove in this article. As we mentioned in the last paragraph, when $\theta_0$ is in $L^2$, the decay (\ref{eq:max-prin-ju}) reduces to the maximum principle (\ref{eq:max-principle}) and no decay rate can be deduced. We address this issue in the following theorems, where we prove that the $L^2$ norm of weak solutions tends to zero but not uniformly, this is, there are solutions with arbitrarily slow decay.

\bt \label{decay-l2} Let $\theta$ be a solution to (\ref{eq:qge-two}) with $\th _0 \in L^2$. Then  

\begin{displaymath}
\lim _{t \to \infty} \Vert \th (t) \Vert _{\eledos} = 0.
\end{displaymath}

\et

\bt \label{decay-slowly} Let $r > 0, \epsilon > 0, T > 0$ be arbitrary. Then, there exists $\th _0$ in $L^2$ with $\| \th _0 \| _{\eledos} = r$ such that if $\theta (t)$ is the solution with  initial data $\th _0$, then 

\begin{displaymath}
\frac{\| \th (T) \| _{\eledos}}{\| \th _0 \| _{\eledos}} \geq 1 - \epsilon.
\end{displaymath}

\et

To prove Theorem \ref{decay-l2} we adapt an argument used in Ogawa, Rajopadhye and Schonbek \cite{o-r-s} to prove decay in the context of the Navier-Stokes equations with slowly varying external forces. It consists in finding estimates for the decay of the norm in the frequency space, studying separately low an high frequencies. The decay of the low frequency part is obtained through generalized energy inequalities, while the Fourier splitting method is used to bound the  decay of the high frequency part. To construct the slowly decaying solutions of Theorem \ref{decay-slowly} we follow the ideas used by Schonbek \cite{sch-3} to prove a similar result for the Navier-Stokes equations. Namely, we construct initial data $\theta_0 ^{\lambda}$ whose $L^2$ norm does not change under an appropiate $\lambda$-scaling, such that it gives rise to a slowly decaying solution to the linear part of (\ref{eq:qge-two}). We then impose extra conditions on  $\theta_0 ^{\lambda}$ to control the term related to the nonlinear part, making it arbitrarily small for small enough values of $\lambda$. 

The next result concerns the decay of the $L^2$ norm of solutions when the initial data is in $L^p \cap \eledos$, with $1 \leq p < 2$.

\bt \label{decay-initial-l2lp}
Let $\th _0 \in L^p \cap \eledos$, where $1 \leq p < 2$. Then, there is a weak solution such that
\begin{displaymath}
\Vert \th (t) \Vert _{\eledos} \leq C (1 + t)^{- \frac{1}{2 \a}(\frac{2}{p} - 1)}.
\end{displaymath}

\et

The proof of Theorem \ref{decay-initial-l2lp} has similarities with  the one  for (\ref{eq:decay-constantin-wu}) in Constantin and Wu \cite{constantin-wu}. We remark that the decay rate we obtain is of the same type as in (\ref{eq:bound-cordoba-cordoba}) and (\ref{eq:max-prin-ju}), where  $\theta_0$ is in $L^1 \cap L^p$, with $ 1 < p < \infty$ and in $L^2 \cap L^p$, with $p \geq 2$, as proved in C\'ordoba and C\'ordoba \cite{cor-cor} and Ju \cite{ju-2} respectively.

The next Theorem is key for establishing decay of the $L^q$ norm of solutions, for large enough $q$. 

\bt \label{kato-like} Let $\Vert \th_0 \Vert _{L^{\exponent}} \leq \kappa$. Then, for $m = \frac{2}{2\a - 1} \leq q < \infty$

\begin{displaymath}
t^{\frac{1}{\a}(\frac{1}{m} - \frac{1}{q})} \th (t) \in BC((0, \infty), L^q).
\end{displaymath}
Moreover, for $\frac{2}{2\a - 1} \leq q < \infty$

\begin{displaymath}
t^{\frac{1}{2 \a} + \frac{1}{\a}(\frac{1}{m} - \frac{1}{q})} \nabla \th (t) \in BC ((0, \infty), L^q).
\end{displaymath}

\et

The proof of Theorem \ref{kato-like} is based on ideas used by Kato \cite{kato} for proving a similar result for the Navier-Stokes equations. Namely, we construct (in an appropiate space) a solution to (\ref{eq:qge-two}) by succesive approximations $\theta _{n+1}$, whose norms are bounded by that of $\theta_1$, $\theta _n$ and $\nabla \theta _n$. This gives rise to a system of recursive inequalities that can be solved if the norm of the data $\theta _0$ is small enough. We can then extract a subsequence converging to a solution with a certain decay rate. This preliminary estimate is then used to obtain the decay rate in Theorem \ref{kato-like}. Note that when $\a = 1$, we recover the rates obtained by Kato \cite{kato}.

\br \, 
A result related to Theorem \ref{kato-like} concerning existence of strong solutions in $L^p$ for small data in $L^q$, where $\frac{2}{2\a - 1} <  q \leq p$ and $\frac{1}{p} + \frac{\a}{q} = \a - \frac{1}{2}$, was proven by Wu \cite{wu-1}. These solutions exist only in an interval $[0,T]$, where the size of the initial data tends to zero when $T$ goes to infinity.  Notice that the special case of $\theta_0$ in  $L^{\exponent}$ is not covered by the hypothesis. 
\er

\br  \label{regularity}
It is well known that solutions to the 3D Navier-Stokes equations, i.e. (\ref{eq:qge-two}) with $\a = 1$, are smooth when  $\theta_0$ is small in $L^2$ and the solution is in $H^1$ (see Heywood \cite{heywood}, Kato \cite{kato} and Serrin \cite{serrin}). For the quasi-geostrophic equation with critical exponent $\a = \frac{1}{2}$, C\'ordoba and C\'ordoba \cite{cor-cor} proved that  when $\theta_0$ is in $H^{\frac{3}{2}}$ and is small in $L^{\infty}$, the solution is in fact classical.  These results suggest that the solution obtained in Theorem \ref{kato-like} might have better regularity than the one obtained. 
\er \\

We now  state the result concerning decay of $L^q$ norms, for large $q$.

\bt \label{theo-decay-2}
Let $ \th _0 \in L^2 \cap L^{\exponent}$. Then there exists $T = T (\th _0)$ such that for $t \geq T$ and $\exponent \leq q < \infty$ 

\begin{displaymath}
\| \theta (t) \| _{L^q} \leq C \, t^{\frac{1}{q}\frac{4a - 3}{\a(2\a - 1)} - 1 + \frac{1}{2 \a}}.
\end{displaymath}

\et

By (\ref{eq:max-prin-ju}), when $ \th _0 \in L^2 \cap L^{\exponent}$, the $L^{\exponent}$ norm of the solution tends to zero. Then, after a (possibly long) time $T = T(\theta _0)$, the solution enters the ball of radius $\kappa$, where $\kappa$ is as in Theorem \ref{kato-like}.  Interpolation between the decays in (\ref{eq:max-prin-ju}) and Theorem \ref{kato-like}, for some $q$ in the appropiate range of values, provides us with a first decay rate. This rate, which is a function of $q$, can then be maximized, leading us to the result in Theorem \ref{theo-decay-2}. A similar idea was used by Carpio \cite{carpio} to obtain analogous results for the Navier-Stokes equations.  

\br \,
After this work was submitted we received preprints of articles by Carrillo and Ferreira \cite{c-f-1}, \cite{c-f-2}, \cite{c-f-3} in which they prove results directly related to the ones obtained here. The proofs by Carrillo and Ferreira are, in general, rather different from ours. In \cite{c-f-2}, they prove Theorem \ref{kato-like} in the particular case $\a = 1$ and $\theta_0 \in \eledos$ but with no restriction on the size of the initial data $\theta_0$. Moreover, they obtain estimates for the decay of all derivatives of $\theta$ in $\eledos \cap L^q$, thus showing that the solution is smooth (see Remark \ref{regularity}). In the forthcoming preprint \cite{c-f-3}, Carrillo and Ferreira extend their results to $\frac{1}{2} < \a \leq 1$ and $\theta_0 \in L^{\exponent}$ and also obtain decays analogous to those of Theorem \ref{theo-decay-2}, but in the more restrictive case of initial data $ \th _0 \in L^1 \cap L^{\exponent}$. 
\er \\

Recently, many articles concerning different aspects of the dissipative quasi-geostrophic equation have been published. Besides the ones we have already referred to, see Berselli \cite{berselli}, Carrillo and Ferreira \cite{c-f-1}, Chae \cite{chae}, Chae and Lee \cite{chae-lee}, Constantin, C\'ordoba and Wu \cite{cons-cor-wu}, Ju \cite{ju-1}, \cite{ju-3},  Schonbek and Schonbek \cite{sch-sch-1}, \cite{sch-sch-2}, Wu \cite{wu-1}, \cite{wu-2}, \cite{wu-3}, \cite{wu-4}, \cite{wu-5}, \cite{wu-6}  and references contained therein.

This article is organized as follows. In Section \ref{preliminaries}, we collect the basic results and estimates we need. In Section \ref{proofs-of theorems-onethree-onfour} we prove Theorems \ref{decay-l2} and \ref{decay-slowly}, in Section \ref{proof-theorem-onetwo} we prove Theorem \ref{decay-initial-l2lp} and finally in Section \ref{proof-of-remaining-theorems} we prove Theorems \ref{kato-like} and \ref{theo-decay-2}.

{\bf Acknowledgments} \, The authors would like to thank Helena Nussenzveig-Lopes for calling their attention to the articles by Carrillo and Ferreira, and Jos\'e Carrillo and Lucas C.F. Ferreira for providing us with copies of their preprints and for very helpful comments and remarks concerning their work.

\section{Preliminaries}
\label{preliminaries}

In this section we collect some essential results and estimates concerning solutions to equation (\ref{eq:qge-two}).

\subsection{Existence and uniqueness of solutions}
\label{basic-facts}

We first state the existence and uniqueness results we assume througout this article.

\bt {\bf (Resnick \cite{resnick})} \, Let $T > 0$ arbitrary. Then, for every $\theta_0 \in \eledos$ and $f \in \eledos([0,T];H^{- \a})$ there exists a weak solution of

\bga
\theta_{t} + (u \cdot \nabla) \theta + (- \Delta)^{\a} \theta = f   \nonumber \\ \th (x,0) = \th_0(x) \nonumber
\eda
such that 

\begin{displaymath}
\theta \in L^{\infty}([0,T];\eledos) \cap \eledos([0,T];H^{\a}).
\end{displaymath}

\et

\bt {\bf (Constantin and Wu  \cite{constantin-wu})} \, Assume that $\a \in (\frac{1}{2},1]$, $T > 0$ and $p$ and $q$ satisfy

\begin{displaymath}
p \geq 1, \quad q > 0, \quad \frac{1}{p} + \frac{\a}{q} = \a - \frac{1}{2}.
\end{displaymath}
Then there is at most one solution $\theta$ of (\ref{eq:qge-two}) with initial value $\theta_0 \in \eledos$ such that

\begin{displaymath}
\theta \in L^{\infty}([0,T];\eledos) \cap \eledos([0,T];H^{\a}), \quad \theta \in L^q([0,T];L^p).
\end{displaymath}
\et

These solutions obey a Maximum Priciple as in (\ref{eq:max-principle}), this is

\begin{displaymath}
\| \theta(t) \| _{L^p} \leq \| \theta_0 \| _{L^p}, \quad t \geq 0
\end{displaymath}
for $ 1 < p \leq \infty$ (see Resnick \cite{resnick}, C\'ordoba and C\'ordoba \cite{cor-cor} and Ju \cite{ju-2} for proofs). Multiplying (\ref{eq:qge-two}) by $\theta$ and integrating in space and time yields 

\be
\label{eqn:int-norm-dif}
\int _s ^t \| \Lambda ^{\a} \theta (\tau) \| _{\eledos} ^2 \, d\tau \leq C, \quad \forall \, s,t > 0. 
\ee

\subsection{Estimates} Let

\be \label{eq:int-form}
\theta (x,t) = K_{\a} (t,x) * \theta _0 - \int _0 ^t K_{\a} (t-s,x) * \nonlinear (s) \, ds
\ee
be the integral form of equation (\ref{eq:qge-two}), where $K_{\a} (x,t)$ is the kernel of the linear part of (\ref{eq:qge-two}), i.e.

\begin{displaymath}
K_{\a}(x,t) = \frac{1}{2 \pi} \int _{\R^2} e^{ix \xi} e^{- |\xi|^{2 \a}t} \, d \xi. 
\end{displaymath}
 
\bp {\bf (Wu \cite{wu-1})} \label{prop-wu}
Let $1 \leq p \leq q \leq \infty$. For any $t > 0$, the operators

\bga
\kernel(t): L^p \to L^q, \quad \kernel (t)\, f = K_\a (t) * f \nonumber \\
\nabla \kernel (t): L^p \to L^q, \quad \nabla \kernel (t) \, f = \nabla K_\a (t) * f \nonumber
\eda
are bounded and 

\bga
\label{eq:norm-ineq-one-lemma}
\Vert \kernel (t)\, f \Vert _{\eleq} \leq C t ^{-\frac{1}{\a}(\frac{1}{p} - \frac{1}{q})} \Vert f \Vert _{L^p} \\
\label{eq:norm-ineq-two-lemma}
\Vert \nabla \kernel (t) \, f \Vert _{\eleq} \leq C t^{-(\frac{1}{2 \a} + \frac{1}{\a}(\frac{1}{p} - \frac{1}{q}))} \Vert f \Vert _{L^p}.
\eda

\ep

The following estimates for the integral term in the right hand side of (\ref{eq:int-form}) are an immediate consequence of Proposition  \ref{prop-wu} and they are key in the proof of  Theorem \ref{kato-like}.

\bl \label{main-lemma-ineq} Let $\eta \leq \mu + \nu < 2$. Then

\be
\label{eq:first-ineq-int-norm}
\Vert \int _0 ^t \kernel (t-s) * \nonlinear (s) \, ds \Vert _{L^{\frac{2}{\eta}}} \leq C \int _0 ^t (t - s)^{-\frac{1}{2a}(\mu + \nu - \eta)} \Vert \theta (s) \Vert _{L^{\frac{2}{\mu}}} \Vert \nabla \theta (s) \Vert _{L^{\frac{2}{\nu}}} ds 
\ee
and

\be
\label{eq:second-ineq-int-norm}
\Vert \int _0 ^t \nabla \kernel (t-s) * \nonlinear (s) \, ds \Vert _{L^{\frac{2}{\eta}}} \leq C \int _0 ^t (t - s)^{-(\frac{1}{2\a}+ \frac{1}{2a}(\mu + \nu - \eta))} \Vert \theta (s) \Vert _{L^{\frac{2}{\mu}}} \Vert \nabla \theta (s) \Vert _{L^{\frac{2}{\nu}}} ds.  
\ee

\el

\pf Use (\ref{eq:norm-ineq-one-lemma}) and (\ref{eq:norm-ineq-two-lemma}) with $q = \frac{2}{\eta}$, $p = \frac{2}{\mu + \nu}$ and

\begin{displaymath}
\Vert \nonlinear \Vert _{L^{\frac{2}{\mu + \nu}}} \leq C \Vert \theta \Vert _{L^{\frac{2}{\mu}}} \Vert \nabla \theta \Vert _{L^{\frac{2}{\mu}}} 
\end{displaymath}
which follows from H\"{o}lder's inequality and boundedness of Riesz transform. $\square$

\bl {\bf (Schonbek and Schonbek \cite{sch-sch-2})} \label{lm:lemma-sch-sch}
Let $\beta, \gamma$ be multi-indices, $|\gamma| < |\beta| + 2 \a  \max (j,1)$, $ j = 0, 1, 2, \cdots $, $1 \leq p \leq \infty$. Then

\begin{displaymath}
\Vert x^{\gamma} D^j _t D^{\beta} \kernel (t) \Vert {L^p} = C t ^{\frac{|\gamma| - |\beta|}{2 \a} - j - \frac{p - 1}{p \a}}
\end{displaymath}
for some constant $C$ depending only on $\alpha, \beta, \gamma, j , p$.
\el

\section{$\eledos$ decay for initial data in $\eledos$}
\label{proofs-of theorems-onethree-onfour}

\subsection{Proof of Theorem \ref{decay-l2}} Let $\theta (t)$ be a solution to (\ref{eq:qge-two}) with $\theta _0 \in \eledos$. For $\phi  = \phi (\xi, t)$

\be \label{eqn:triangle-norm}
\| \hattheta (t) \| _{\eledos} ^2 \leq 2 \left( \|\phi (t) \hattheta (t) \| _{\eledos} ^2 + \| (1 - \phi (t)) \hattheta (t) \| _{\eledos} ^2 \right).
\ee
We call the terms $\|\phi (t) \hattheta (t) \| _{\eledos} ^2$ and $\| (1 - \phi (t)) \hattheta (t) \| _{\eledos} ^2$ the {\em low and high frequency parts of the energy} respectively. In Propositions  \ref{prop-low-freq} and \ref{prop-high-freq} and Corollary \ref{cor-low-freq}, we obtain estimates, for an appropiate class of functions $\phi$, that allow us to prove that the low and high frequency parts of the energy tend to zero. These estimates are of similar character  to the ones that Ogawa, Rajophadye and Schonbek obtained for the Navier-Stokes equations in \cite{o-r-s}. 

\subsubsection{Energy estimates} We first establish some preliminary estimates which will be needed in the proof of  Theorem \ref{decay-l2}.

\bp \label{prop-low-freq}
Let $\psi \in C^1((0,\infty), C^1 \cap \eledos)$. Then for $0 < s <t$

\bga
\| \hattheta \hatpsi (t) \| _{\eledos} ^2 &  \leq & \| \hattheta \hatpsi (s) \| _{\eledos} ^2 + 2 \int _s ^t | \langle \hatpsi ' \hattheta (\tau), \hatpsi \hattheta (\tau) \rangle - \| |\xi| ^{\a} \hattheta \hatpsi (\tau) \| _{\eledos} ^2 | \, d \tau \nonumber \\ & + & 2 \int _s ^t | \langle \xi \cdot \widehat{u \theta} (\tau), \hatpsi ^2 \hattheta (\tau)\rangle | \, d \tau. \nonumber
\eda 
\ep

\pf Let $\theta (t)$ be a smooth solution to (\ref{eq:qge-two}). Taking the Fourier transform, multiplying by $\hatpsi ^2 \hattheta$ and integrating by parts we obtain the formal estimate

\bga
\frac{d}{dt} \| \hatpsi \hattheta \| _{\eledos} ^2 & = & 2 \left( \langle \hatpsi ' \theta (t), \hatpsi \hattheta (t) \rangle - \| |\xi|^{\a} \hatpsi \hattheta (t) \| _{\eledos} ^2 \right) \nonumber \\ & - & 2 \langle \widehat{\nonlinear} (t), \hatpsi ^2 \hattheta (t) \rangle. \nonumber
\eda
Integrating between $s$ and $t$ yields

\bga
\| \hattheta \hatpsi (t) \| _{\eledos} ^2 &  \leq & \| \hattheta \hatpsi (s) \| _{\eledos} ^2 + 2 \int _s ^t | \langle \hatpsi ' \hattheta (\tau), \hatpsi \hattheta (\tau) \rangle - \| |\xi| ^{\a} \hattheta \hatpsi (\tau) \| _{\eledos} ^2 | \, d \tau \nonumber \\ & + & 2 \int _s ^t | \langle \xi \cdot \widehat{u \theta} (\tau), \hatpsi ^2 \hattheta (\tau)\rangle | \, d \tau. \nonumber
\eda 
As before, the retarded mollifiers method allows us to extend this estimate to weak solutions.  For full details see Ogawa, Rajopadhye and Schonbek \cite{o-r-s}. $\square$

\bc \label{cor-low-freq}
Let $\phi \in C^1((0, \infty), L^2)$. Then for $0 < s < t$

\begin{displaymath}
\| \hattheta \phi (t) \| _{\eledos} ^2 \leq \| \hattheta (s) e^{-|\xi|^{2\a}(t-s)} \phi (t) \| _{\eledos} ^2 + 2 \int _s ^t | \langle \xi \cdot \widehat{u \theta}, e^{-2|\xi|^{2\a}(t-\tau)} \phi ^2 (t) \hattheta (\tau) \rangle| \, d \tau.
\end{displaymath}
\ec

\pf Take $\hatpsi _{\eta} (\tau) = e^{-|\xi|^{2 \a}(t + \eta - \tau)} \phi(\xi,t)$ for $\eta > 0$. Then

\begin{displaymath}
\langle \hatpsi _{\eta} ' \hattheta (\tau), \hatpsi _{\eta} \hattheta (\tau) \rangle = \langle |\xi|^{2 \a} \hatpsi _{\eta} \hattheta (\tau),  \hatpsi _{\eta} \hattheta (\tau) \rangle = \| |\xi| ^{\a} \hatpsi _{\eta} \hattheta (\tau) \| _{\eledos} ^2
\end{displaymath}
so the integrand in the second term in Proposition \ref{prop-low-freq} vanishes. Taking limit as $\eta \to 0$ we see that $\hatpsi (t) = \phi(\xi, t)$ and $\hatpsi (s) = e^{-|\xi|^{2 \a}(t - s)} \phi (\xi,t)$, so

\begin{displaymath}
\| \hattheta \phi (t) \| _{\eledos} ^2 \leq \| \hattheta (s) e^{-|\xi|^{2\a}(t-s)} \phi (t) \| _{\eledos} ^2 + 2 \int _s ^t | \langle \xi \cdot \widehat{u \theta}, e^{-2|\xi|^{2\a}(t-\tau)} \phi ^2 (t) \hattheta (\tau) \rangle| \, d \tau.
\end{displaymath}
as we wanted to prove. $\square$

\bp \label{prop-high-freq}
Let $E \in C^1 ((0, \infty), \R)$ and $\psi \in C^1((0,\infty), L^{\infty})$ such that $1 - \psi^2 \in L^{\infty} ((0,\infty),L^{\infty})$ and $\nabla \Ff ^{-1} (1 - \psi^2) \in L^{\infty} ((0,\infty),\eledos)$. Then

\bga
E(t) \| \psi \hattheta (t) \| _{\eledos} ^2 & \leq & E(s) \| \psi \hattheta (s) \| _{\eledos} ^2 + \int _s ^t E'(\tau) \| \psi \hattheta (\tau) \| _{\eledos} ^2 d \tau \nonumber \\ & + & 2 \int _s ^t E(\tau) | \langle \psi ' \hattheta (\tau), \psi \hattheta (\tau) \rangle - \| |\xi| ^{\a} \psi \hattheta (\tau )\| _{\eledos} ^2 | \, d \tau \nonumber \\ & + & 2 \int _s ^t E(\tau) |\langle \widehat{\nonlinear} (\tau), (1 - \psi^2 (\tau)) \hattheta (\tau) \rangle| \, d \tau. \nonumber
\eda
\ep 

\pf We prove the estimate first for smooth solutions. As in Proposition \ref{prop-low-freq}, we take the Fourier transform of (\ref{eq:qge-two}) and multiply it by $E \psi ^2 \hattheta$. Integrating by parts and then between $s$ and $t$ we obtain the formal estimate

\bga
E(t) \| \psi \hattheta (t) \| _{\eledos} ^2 & \leq & E(s) \| \psi \hattheta (s) \| _{\eledos} ^2 + \int _s ^t E'(\tau) \| \psi \hattheta (\tau) \| _{\eledos} ^2 d \tau \nonumber \\ & + & 2 \int _s ^t E(\tau) | \langle \psi ' \hattheta (\tau), \psi \hattheta (\tau) \rangle - \| |\xi| ^{\a} \psi \hattheta (\tau )\| _{\eledos} ^2 | \, d \tau \nonumber \\ & + & 2 \int _s ^t E(\tau) |\langle \widehat{\nonlinear} (\tau), (1 - \psi^2 (\tau)) \hattheta (\tau) \rangle| \, d \tau. \nonumber
\eda
Here we used that $\langle (\nonlinear, \theta \rangle = 0$. When using the retarded mollifiers method, the conditions $1 - \psi^2 \in L^{\infty} ((0,\infty),L^{\infty})$ and $\nabla \Ff ^{-1} (1 - \psi^2) \in L^{\infty} ((0,\infty),\eledos)$ will guarantee the weak convergence of the nonlinear term. For full details see Ogawa, Rajopadhye and Schonbek \cite{o-r-s}. $\square$ 

\subsubsection{Proof of Theorem \ref{decay-l2}} We first prove the following easy estimate.

\begin{Lemma} \label{lm:decay-fm}
Let $m > 0$ and 

\begin{displaymath}
f_m (t) = \int _{|\xi| > 1} |\xi| ^{2\a} e^{- m |\xi|^{2 \a}t} \, d \xi. 
\end{displaymath}
Then $\lim _{t \to \infty} f_m (t) = 0$. 
\end{Lemma}

\pf From the inequality

\begin{displaymath}
|\xi| ^{2\a} e^{- m |\xi|^{2 \a}t} \leq C \ \frac{e^{-m |\xi|^{2 \a}t/2}}{mt}
\end{displaymath}
it follows that

\bga \label{eqn:bound-lastterm-fm}
f_m (t) & = & \int _{|\xi| > 1} |\xi| ^{2\a} e^{- m |\xi|^{2 \a}t} \, d \xi \leq C \int _{|\xi| > 1} \frac{e^{-m |\xi|^{2 \a}t/2}}{mt} \, d \xi \nonumber \\ & \leq & C \int _{|\xi| > 1} \frac{e^{-m |\xi| t/2}}{mt} \, d \xi = C \int _1 ^{\infty} \frac{r e^{-mtr/2}}{mt} \, dr \nonumber \\ & = & C \frac{e^{-mtr/2}}{(mt)^2} \left( 1 + \frac{2}{mt} \right) \leq \frac{C}{(mt)^2}.
\eda
Thus, $\lim _{t \to \infty} f_m (t) = 0$. $\square$

We choose $\phi(\xi,t) = e^{-|\xi|^{2\a}t}$. Note that $\phi$ is the kernel of the solution to the Fourier transform of (\ref{eq:qge-two}). \\

{\em Low frequency energy decay.} Using Corollary \ref{cor-low-freq} with $\phi$ as defined above we obtain

\be \label{eqn:low-freq}
\| \hattheta \phi (t) \| _{\eledos} ^2 \leq \| \hattheta (s) \phi(t-s) \phi(t) \| _{\eledos} ^2 + 2 \int _s ^t | \langle \xi \cdot \widehat{u \theta}, \phi^2 (t - \tau) \phi ^2 (t) \hattheta (\tau) \rangle.
\ee
A standard application of the Dominated Convergence Theorem proves that the first term in the right hand side of (\ref{eqn:low-freq}) tends to zero when $t$ goes to infinity. Now

\bga \label{eqn:ineq-inner-prod}
\langle \xi \cdot \widehat{u \theta} (\tau), \phi ^2 (t-\tau) \phi ^2 (t) \hattheta (\tau) \rangle & = & \langle \xi^{\a} \cdot \widehat{u \theta} (\tau), |\xi| ^{1 - \a} \phi ^2(t-\tau) \phi ^2 (t) \hattheta (\tau) \rangle \nonumber \\ & \leq  & \| |\xi| ^{\a} \widehat{u \theta} (\tau) \| _{L^{\infty}} \| |\xi| ^{1 -\a} \phi ^2(t-\tau) \phi ^2 (t) \hattheta (\tau) \| _{L^1}.
\eda
As $\widehat{u \theta} (\xi) = \hat{u} * \hattheta (\xi)$, then 

\bga \label{eqn:linfinity-norm} 
\| |\xi| ^{\a} \widehat{u \theta} (\tau) \| _{L^{\infty}} & = & \sup _{\xi \in \R^2} \left| \int _{\R^2} |\xi|^{\a} \hat{u} (\xi - \eta) \hattheta (\eta) \, d \eta \right| \nonumber \\ & \leq & C \left( \sup _{\xi \in \R^2} \left| \int _{\R^2} |\xi - \eta|^{\a} \hat{u} (\xi - \eta) \hattheta (\eta) \, d \eta \right| + \sup _{\xi \in \R^2} \left| \int _{\R^2} |\eta|^{\a} \hat{u} (\xi - \eta) \hattheta (\eta) \, d \eta \right| \right) \nonumber \\ & \leq & C \left( \| \widehat{(\Lambda ^{\a})u \, \theta }\| _{L^{\infty}} + \| \widehat{u \, \Lambda ^{a} \theta} \| _{L^{\infty}} \right) \leq C \left( \| (\Lambda ^{\a}u) \theta \| _{L^1} + \|u \, \Lambda ^{\a} \theta \| _{L^1} \right) \nonumber \\ & \leq &  C \| \Lambda ^{\a} \theta \| _{\eledos}. \nonumber
\eda
Now 

\be \label{eqn:bound-lone-norm}
\| |\xi| ^{1 -\a} \phi ^2 (t-\tau) \phi ^2 (t) \hattheta (\tau) \| _{L^1}  \leq  C \| |\xi| ^{1 -2\a} \phi ^2 (t-\tau) \phi ^2 (t) \| _{\eledos} \| |\xi| ^{\a} \hattheta (\tau) \| _{\eledos}.
\ee
and

\bga \label{eqn:bound-ltwo-norm}
\| |\xi| ^{1 -2\a} \phi ^2 (t-\tau) \phi ^2 (t) \| ^2 _{\eledos} & = & \int _{\R ^2} |\xi| ^{2-4\a} e^{-4 | \xi |^{2 \a}(2t - \tau)} \, d \xi \nonumber \\ & \leq & \int _{|\xi| \leq 1} | \xi | ^{2-4\a} \, d \xi + f_4 (2t - \tau).
\eda
As $ \frac{1}{2} < \a \leq 1$, the first term in the right hand side of (\ref{eqn:bound-ltwo-norm}) is integrable. By Lemma \ref{lm:decay-fm} and (\ref{eqn:bound-lastterm-fm}) we see that

\be \label{eqn:f4m-bounded}
f_4(2t - \tau) \leq \frac{C}{4(2t - \tau)^2} \leq C
\ee
for $t$ large enough. Then using (\ref{eqn:ineq-inner-prod}), (\ref{eqn:linfinity-norm}), (\ref{eqn:bound-lone-norm}), (\ref{eqn:bound-ltwo-norm}) and (\ref{eqn:f4m-bounded}) we obtain that

\bga
2 \int _s ^t | \langle \xi \cdot \widehat{u \theta}, \phi ^2 (t-\tau) \phi ^2 (t) \hattheta (\tau) \rangle \, d \tau & \leq & \int _s ^t \| |\xi| ^{\a} \widehat{u \theta} (\tau) \| _{L^{\infty}} \| |\xi| ^{1 -\a} \phi ^2 (t-\tau) \phi ^2 (t) \hattheta (\tau) \| _{L^1} \, d \tau \nonumber \\ & \leq & 2 \int _s ^t |\xi| ^{\a} \widehat{u \theta} (\tau) \| _{L^{\infty}} \| |\xi| ^{1 -2\a} \phi ^2 (t-\tau) \phi ^2 (t) \hattheta (\tau) \| _{\eledos} \| |\xi| ^{\a} \hattheta (\tau) \| _{\eledos} \, d \tau  \nonumber \\ & \leq & C \int _s ^t \| \Lambda  ^{\a} \theta (\tau) \| _{\eledos} ^2 \, d \tau. \nonumber
\eda
Taking limits as $s$ and $t$ go to infinity, (\ref{eqn:int-norm-dif}) implies that the low frequency part of the energy goes to zero. \\

{\em High energy frequency decay.} Let $\psi(\xi,t) = 1 - \phi(\xi,t)$. As 
\begin{displaymath}
1 - \psi ^2(\xi,t) = 2 \phi(\xi,t) - \phi^2 (\xi,t) = 2 e^{-|\xi|^{2\a}t} - e^{-2|\xi|^{2\a}t}
\end{displaymath}
 decays exponentially fast, we can apply Proposition \ref{prop-high-freq}. After rearranging terms, we obtain 

\begin{displaymath}
\| (1 - \phi(t)) \hattheta (t) \| _{\eledos} ^2 = \| \psi (t) \hattheta (t) \| _{\eledos} ^2 \leq I + II + III + IV
\end{displaymath}
where

\begin{displaymath}
I = \frac{E(s)}{E(t)} \| \psi \hattheta (s) \| _{\eledos} ^2
\end{displaymath}

\begin{displaymath}
II = \frac{1}{E(t)} \int _s ^t \left( E'(\tau) \| \psi  \hattheta (\tau) \| _{\eledos} ^2 - 2 E(\tau) \| |\xi| ^{\a} \psi \hattheta (\tau) \| _{\eledos} ^2 \right) \, d \tau
\end{displaymath}

\begin{displaymath}
III = \frac{2}{E(t)} \int_s ^t E(\tau) \langle \psi ' \hattheta (\tau), \psi (\tau) \hattheta (\tau) \rangle \, d \tau
\end{displaymath}

\begin{displaymath}
IV = \frac{2}{E(t)} \int _s ^t E(\tau) | \langle \widehat{\nonlinear} (\tau), (1 - \psi^2 (\tau)) \hattheta (\tau) \rangle | \, d \tau.
\end{displaymath}
We choose $E(t) = \left( 1 + t \right) ^k$, where $k > 2$. \\

\subparagraph{Term I} Since $|\psi| \leq C$ and $\theta \in L^2$

\begin{displaymath}
I = \left( \frac{1 + s}{1 + t} \right)^k \| \psi (s) \hattheta (s) \| _{\eledos} ^2 \leq C \left( \frac{1 + s}{1 + t} \right)^k.
\end{displaymath}
Thus, 

\begin{displaymath}
\lim _{t \to \infty} I(t) = 0.
\end{displaymath}

\subparagraph{Term II} We use the Fourier splitting method. Let 
\begin{displaymath}
B(t) = \{ \xi \in \R^2: |\xi| \leq G(t) \}
\end{displaymath}
where $G$ is to be determined below. Then 

\bga \label{eqn:long-ineq}
E'(\tau) \| \psi (\tau) \hattheta (\tau) \| _{\eledos} ^2 & - & 2 E(\tau) \| |\xi| ^{\a} \psi (\tau) \hattheta (\tau) \| _{\eledos} ^2  \nonumber \\ & = & E'(\tau) \int _{\R^2 \smallsetminus B(t)} | \psi \hattheta (\tau) | ^2 \, d \xi - 2 E(\tau)  \int _{\R^2 \smallsetminus B(t)} |\xi| ^{2\a} |\psi \hattheta (\tau)| ^2 \, d \xi \nonumber \\ & + & E'(\tau) \int _{B(t)} | \psi \hattheta (\tau) | ^2 \, d \xi - 2 E(\tau)  \int _{B(t)} |\xi| ^{2\a} |\psi \hattheta (\tau)| ^2 \, d \xi \nonumber \\ & \leq & \left( E'(\tau) - 2 E(\tau) G^{2 \a} (\tau) \right) \int _{\R^2 \smallsetminus B(t)} | \psi \hattheta (\tau) | ^2 \, d \xi \nonumber \\ & + & E'(\tau) \int _{B(t)} | \psi \hattheta (\tau) | ^2 \, d \xi - 2 E(\tau)  \int _{B(t)} |\xi| ^{2\a} |\psi \hattheta (\tau)| ^2 \, d \xi.
\eda
Choosing $G(t) = \left( \frac{k}{2(1+t)} \right) ^{\frac{1}{2 \a}}$, we see that $E'(\tau) - 2 E(\tau) G^{2 \a} (\tau) = 0$ so the first term in the right hand side of (\ref{eqn:long-ineq}) vanishes. As the last term in (\ref{eqn:long-ineq}) is negative, it can be dropped, hence

\begin{displaymath}
II \leq \frac{k}{(1 + t)^k} \int _s ^t (1 + \tau) ^{k -1} \left( \int _{B(t)}  |\psi \hattheta (\tau)| ^2 \, d \xi \right) \, d \tau.
\end{displaymath}
As $\psi (\xi,t) = 1 - e^{-|\xi|^{2\a}t}$, then $|\psi| \leq |\xi|^{2 \a}$ for $|\xi| \leq 1$. Then

\begin{displaymath}
\int _{B(t)}  |\psi \hattheta (\tau)| ^2 \, d \xi \leq \int _{B(t)} |\xi|^{4 \a}  |\hattheta (\tau)| ^2 \, d \xi \leq C G^{4 \a} (t) = \frac{C}{(1 + t)^2}.
\end{displaymath}
Then

\begin{displaymath}
II \leq \frac{k}{(1 + t)^k} \int _s ^t (1 + \tau) ^{k - 3} \, d \tau \leq \frac{C}{(1 + t)^2},
\end{displaymath}
so

\begin{displaymath}
\lim _{t \to \infty} II (t) = 0.
\end{displaymath}

\subparagraph{Term III} As $\psi (\xi,t) = 1 - e^{-|\xi|^{2\a}t}$, then $\psi' = \frac{\partial \psi}{\partial t} = |\xi|^{2 \a} e^{-|\xi|^{2\a}t} = |\xi|^{2 \a} \phi (\xi,t)$. As $E(t)$ is an increasing function

\bga \label{eqn:bound-three}
III & = & \frac{2}{E(t)} \int_s ^t E(\tau) \langle |\xi|^{2 \a} \phi(\tau)  \hattheta (\tau), (1 - \phi(\tau))  \hattheta (\tau) \rangle \, d \tau   \nonumber \\ & \leq & 2 \int_s ^t  \langle |\xi|^{\a}  \hattheta (\tau), |\xi|^{\a} \hattheta (\tau) \rangle \, d \tau  \nonumber \\ & \leq & 2 \int _s ^t \| \Lambda ^{\a} \hattheta (\tau) \| _{\eledos} ^2 \, d \tau.
\eda
Taking limits  when $t$ and $s$ go to infinity we obtain

\begin{displaymath}
\lim _{t \to \infty} III (t) = 0.
\end{displaymath}

\subparagraph{Term IV} Let $\hat{\omega} (\xi, t) = 1 - \psi^2(\xi,t)$. Then

\bga \label{eqn:first-bound-term4}
| \langle \widehat{\nonlinear} (\tau), \hat{\omega} (\tau) \hattheta (\tau) \rangle | & \leq  & \langle |\xi|^{\a} |\widehat{u \theta} (\tau)|, |\xi|^{1 -a} |\hat{\omega} \hattheta( \tau)| \rangle \nonumber \\ & \leq & \| |\xi| ^{\a} \widehat{u \theta} (\tau) \| _{L^{\infty}} \| |\xi| ^{1 -\a} \hat{\omega} \hattheta (\tau) \| _{L^1}.
\eda

As in (\ref{eqn:bound-lone-norm})

\be \label{eqn:second-bound-term4}
\| |\xi| ^{1 -\a} \hat{\omega} (\tau) \hattheta (\tau) \| _{L^1}  \leq  C \| |\xi| ^{1 -2\a} \hat{\omega} (\tau) \| _{\eledos} \| |\xi| ^{\a} \hattheta (\tau) \| _{\eledos}.
\ee
We notice that as $2\a - 4 < 0$, then

\bga \label{eqn:third-bound-term4}
\Vert |\xi| ^{1 -2\a} \hat{\omega} (\tau) \| _{\eledos} ^2 & = & \int _{\R ^2} |\xi| ^{2 - 4 \a} |\hat{\omega}|^2 \, d \xi \nonumber \\ & \leq & \int _{|\xi| \leq 1} |\xi| ^{2 - 4 \a} \, d \xi + \int_{|\xi| \geq 1} |\hat{\omega}|^2 d \xi \leq C.
\eda
 Then by (\ref{eqn:first-bound-term4}), (\ref{eqn:second-bound-term4}) and (\ref{eqn:third-bound-term4}) 

\bga
IV & \leq & 2 \int _s ^t |\langle \widehat{\nonlinear} (\tau), \hat{\omega} (\tau) \hattheta (\tau) \rangle| \, d \tau \nonumber \\ & \leq & 2 \int _s ^t \| |\xi| ^{\a} \widehat{u \theta} (\tau) \| _{L^{\infty}} \| |\xi| ^{1 -2\a} \hat{\omega} (\tau) \| _{\eledos} \| \Lambda  ^{\a} \theta (\tau) \| _{\eledos} \, d \tau \nonumber \\ & \leq & 2 \int _s ^t \| \Lambda  ^{\a} \theta (\tau) \| ^2 _{\eledos} \, d \tau. \nonumber 
\eda
As before, letting $s$ and $t$  go to infinity we obtain

\begin{displaymath}
\lim _{t \to \infty} IV(t) = 0.
\end{displaymath} 
Thus, the high frequency part of the energy goes to zero, which concludes the proof of Theorem \ref{decay-l2}. $\square$

\subsection{Proof of Theorem \ref{decay-slowly}} We briefly describe the idea of the proof. In order to make the decay of a solution to (\ref{eq:qge-two}) arbitrarily slow, we will construct a set of initial data $\{ \theta _0 ^{\lambda} \} _{\lambda > 0}$ in $L^2$  such that $\Vert \theta _0 ^{\lambda} \Vert _{\eledos} = \Vert \theta _0  \Vert _{\eledos}$. The mild solution to (\ref{eq:qge-two}) with initial data $\theta _0 ^{\lambda}$  

\be \label{eq:scaled-equation}
\theta ^{\lambda} (x,T) =  \kernel (T) * \theta _0 ^{\lambda} (x) - \int _0 ^T \kernel (T-s) * (u ^{\lambda} \cdot \nabla) \theta ^{\lambda} (s) \, ds
\ee
has the following property: given $T > 0$, we can find  $\lambda$ sufficiently close to zero, so that the $L^2$ norm of the first term of the right hand side of (\ref{eq:scaled-equation}) stays arbitrarily close to that of $\theta _0$. For this to hold, $\theta _0 ^{\lambda}$ must be such that: $a)$ the $L^2$ norm of $\theta _0 ^{\lambda}$ is invariant under the scaling; $b)$ $\theta _0 ^{\lambda}$ gives rise to a self-similar solution to the linear part of (\ref{eq:qge-two}); and $c)$ $\theta _0$ is in $L^p \cap L^q$, for adecuate $p$ and $q$, so that he integral term will be sufficientyl small. We remark that as a result of our choice, the $L^p$ and $L^q$ norms of $\theta _0 ^{\lambda}$ will not be invariant under scaling. 

We proceed to the proof now. For $\theta _0$ in $L^2$, it is easy to see that $\theta _0 ^{\lambda} (x) = \lambda \theta_0 (\lambda x)$ is such that 

\begin{displaymath}
\| \theta _0 ^{\lambda} \| _{\eledos} = \| \theta _0 \| _{\eledos}, \quad \lambda > 0.
\end{displaymath}
Then, for these $\theta _0 ^{\lambda}$, condition $a)$ holds. Now let $\theta _0$ be such that $\theta _0 ^{\lambda}$ gives rise to a self-similar solution $\Theta ^{\lambda}$ to the linear part of (\ref{eq:qge-two}), this is

\begin{displaymath}
\Theta ^{\lambda} (x,t) = \lambda \Theta (\lambda x, \lambda ^{2 \a} t)
\end{displaymath}
is a solution to 

\bga
\Theta _t + (- \Delta)^{\a} \Theta = 0 \nonumber \\ \Theta ^{\lambda} _0 (x) = \theta _0 ^{\lambda} (x). \nonumber 
\eda
By uniqueness of the solution to the linear part, we have  $\Theta ^{\lambda} (x,t) = K_{\a}(t) * \Theta ^{\lambda} _0 (x) = K_{\a}(t) *\theta _0 ^{\lambda} (x)$, thus

\bga
\| \Theta ^{\lambda} (t) \| _{\eledos} & = & \int _{\R ^2} | \Theta ^{\lambda} (x,t)| ^2 \, dx = \lambda ^2  \int _{\R ^2} | \Theta  (\lambda x,\lambda ^{2\a} t)| ^2 \, dx \nonumber \\ & = & \int _{\R ^2} | \Theta  (y,\lambda ^{2 \a} t)| ^2 \, dy  = \int _{\R ^2}  e ^{- |\xi|^{2 \a} \lambda ^{2\a} t} |\hat{\theta} _0 (\xi)| \, d \xi. \nonumber
\eda
As a result of this, given $T > 0$

\be \label{eq:limit-linear-term}
\lim _{\lambda \to 0} \frac{\| \widehat{\Theta ^{\lambda}} (T) \| _{\eledos} ^2 }{\| \widehat{\Theta  _0} \| _{\eledos} ^2} = \lim _{\lambda \to 0} \frac{\int _{\R ^2}  e ^{- |\xi|^{2 \a} \lambda ^{2\a} T} |\hat{\theta} _0 (\xi)|^2  \, d \xi}{\int _{\R^2} |\widehat{\theta _0} (\xi)|^2 \, d \xi} = 1.
\ee
This shows that choosing $\lambda$ small enough, we can make the ratio of the norms arbitrarily close to $1$ for large enough $t$.

We now address the integral term in (\ref{eq:scaled-equation}). We first notice that

\bga
\| \kernel (t-s) * (u ^{\lambda} \cdot \nabla) \theta ^{\lambda} (s) \| _{\eledos} & = & \| \nabla \kernel (t-s) * (u ^{\lambda} \theta ^{\lambda}) (s) \| _{\eledos} \nonumber \\ & \leq & \| \nabla \kernel (t-s) \| _{L^1} \| (u ^{\lambda} \theta ^{\lambda}) (s)\| _{\eledos} \nonumber \\ & \leq & C (t-s) ^{-\frac{1}{2 \a}} \| u ^{\lambda} (s) \|_{L^p} \| \theta ^{\lambda} (s) \| _{L^q}  \nonumber \\ & \leq &  C (t-s) ^{-\frac{1}{2 \a}} \| \theta ^{\lambda} (s) \|_{L^p} \| \theta ^{\lambda} (s) \| _{L^q} \nonumber
\eda
where we have used Lemma  \ref{lm:lemma-sch-sch} with $\gamma = 0, p = 1, \beta = 1, j = 0$,  H\"older's inequality with $\frac{1}{2} = \frac{1}{p} + \frac{1}{q}$ and boundedness of the Riesz transform. By the Maximum Principle (\ref{eq:max-principle}), $\| \theta ^{\lambda} (s) \| _{L^m} \leq \| \theta ^{\lambda} _0\| _{L^m}$ and as 

\be \label{eq:norm-scaled-function}
\| \theta ^{\lambda} _0 \| _{L^m} = \lambda ^{1 - \frac{2}{m}} \| \theta _0\| _{L^m}
\ee
then

\bga
\| \kernel (t-s) * (u ^{\lambda} \cdot \nabla) \theta ^{\lambda} (s) \| _{\eledos} & \leq & C (t-s) ^{-\frac{1}{2 \a}} \lambda ^{2 - (\frac{2}{p} + \frac{2}{q})} \| \theta _0\| _{L^p} \| \theta _0\| _{L^q} \nonumber \\ & \leq & C (t-s) ^{-\frac{1}{2 \a}} \lambda \| \theta _0\| _{L^p} \| \theta _0\| _{L^q}.
\eda
We remark that by (\ref{eq:norm-scaled-function}), the $L^m$ norm of  $\theta ^{\lambda} _0$ is invariant only when $m = 2$. Choosing $\theta _0$ in $L^p \cap L^q$ (condition $c)$) we obtain

\be \label{eq:bound-integral-lambda}
\int _0 ^T \| \kernel (T-s) * (u ^{\lambda} \cdot \nabla) \theta ^{\lambda} (s) \| _{\eledos} \, ds  \leq  C T ^{1 - \frac{1}{2 \a}} \lambda \| \theta _0\| _{L^p} \| \theta _0\| _{L^q}.
\ee
So given $\epsilon > 0$ and $T > 0$, we can choose $\lambda > 0$ such that by (\ref{eq:limit-linear-term})

\begin{displaymath}
\frac{\| \kernel (T) * \theta _0 ^{\lambda} \| _{\eledos}}{\| \theta _0 ^{\lambda}\| _{\eledos}} \geq 1 - \frac{\epsilon}{2}
\end{displaymath}
and by (\ref{eq:bound-integral-lambda})

\begin{displaymath}
\frac{\int _0 ^T \| \kernel (T-s) * (u ^{\lambda} \cdot \nabla) \theta ^{\lambda} (s) \| _{\eledos} \, ds}{\| \theta _0 ^{\lambda}\| _{\eledos}} \leq \frac{\epsilon}{2}.
\end{displaymath}
Then

\begin{displaymath}
\frac{\| \th ^{\lambda} (T) \| _{\eledos}}{\| \th ^{\lambda} _0 \| _{\eledos}} \geq 1 - \epsilon.
\end{displaymath}
This proves our result.  $\square$

\section{$L^2$ decay for initial data in $L^p \cap \eledos, 1 \leq p < 2$}
\label{proof-theorem-onetwo}

To prove Theorem \ref{decay-initial-l2lp}, we follow a modified version of the Fourier splitting method, see Constantin and Wu \cite{constantin-wu}. Similar ideas in the context of the 2D Navier-Stokes equation can be found in Zhang \cite{zhang}. In order to compute the actual decay rate of the $L^2$ norm, we need a preliminary estimate, proven in Lemma \ref{lm:preliminary-estimate}, which we then use to establish the right decay. In both proofs we first obtain formal estimates for smooth solutions through the Fourier splitting method and we  then use the method of retarded mollifiers of Cafarelli, Kohn and Nirenberg \cite{ckn} to extend them to weak solutions.

The following auxiliary Lemmas will be necessary in the sequel.

\begin{Lemma}  Let $h \in L^p$, $1 \leq p < 2$ and let $S(t) = \{ \xi \in \R^2: |\xi| \leq g(t) ^{-\frac{1}{2\a}}\}$, for a continuous function $g:\R ^+ \to \R ^+$. Then

\begin{displaymath}
\int _{S(t)} |\hat{h}| ^2 d \xi \leq C g(t) ^{-\frac{1}{\a}(\frac{2}{p}  -1)}.
\end{displaymath}
\label{lm:lemma-1}
\end{Lemma}

\pf By Cauchy-Schwarz

\begin{displaymath}
\int _{S(t)} |\hat{h}| ^2 d \xi \leq \left( \int _{S(t)} |\hat{h}| ^{2r} d \xi \right) ^{\frac{1}{r}} \left( \int _{S(t)}d \xi \right) ^{\frac{1}{s}}
\end{displaymath}
where $\frac{1}{r} + \frac{1}{s} = 1$. Setting $2r = q$, we obtain $\frac{1}{r} = \frac{2}{q}$ and $\frac{1}{s} = \frac{2}{p} - 1$. By the Riesz-Thorin Interpolation Theorem, $\Ff: L^p \ \to L^q$ is bounded for $p \in [1,2]$ and $\frac{1}{p} + \frac{1}{q} = 1$. As $h$ is in $L^p$, then $\| \hat{h} \| _{\eleq} \leq \Vert h \Vert _{L^p}$ and as a result of this

\bga
\int _{S(t)} |\hat{h}| ^2 d \xi & \leq & C \left( \int _{S(t)}d \xi \right) ^{\frac{2}{p}-1} = C \, (Vol \, S(t))^{\frac{2}{p}-1} \nonumber \\ & = & C \, r(t) ^{2(\frac{2}{p}-1)}  = C g(t) ^{-\frac{1}{\a}(\frac{2}{p}  -1)}. \quad \square \nonumber
\eda

\begin{Lemma} Let $\theta$ be a solution to (\ref{eq:qge-two}). Then, 

\begin{displaymath}
|\widehat{\nonlinear} (\xi)| \leq C |\xi| \| \theta \| _{\eledos} ^2.
\end{displaymath}
\label{lm:lemma-2}
\end{Lemma}

\pf As $\widehat{\nonlinear} (\xi) = \widehat{\nabla \cdot u \theta} (\xi) =   \xi \cdot \widehat{u \theta} (\xi)$, boundedness of the Fourier transform and of the Riesz transform imply

\begin{displaymath}
|\widehat{\nonlinear} (\xi)| = | \xi | | \widehat{u \theta} (\xi)| \leq |\xi| \Vert \widehat{u \theta} \Vert _{L^{\infty}} \leq C | \xi | \Vert u \theta \Vert _{L^1} \leq C | \xi | \Vert \theta \Vert _{\eledos} ^2. \quad \square
\end{displaymath}

In the next Lemma we establish the preliminary decay rate.

\begin{Lemma} \label{lm:preliminary-estimate}
Let $\theta$ be a solution to (\ref{eq:qge-two}) with initial data $\theta_0$ in $L^p \cap \eledos, 1 \leq p < 2$. Then

\begin{displaymath}
\int _{\R^2} | \hat{\theta} | ^2 d \xi \leq C \ln (e + t)^{-(1 + \frac{1}{\alpha})}.
\end{displaymath}

\end{Lemma}

\pf The first part of the proof consists of a formal argument that proves the expected decay for smooth solutions. At the end of the proof we sketch how to make the argument rigurous. We use the Fourier splitting method, taking

\begin{displaymath}
B(t) = \{ \xi: |\xi| \leq g ^{-\frac{1}{2\a}} (t) \}
\end{displaymath}
where $g(t) = (\frac{1}{2} + \frac{1}{2 \a}) [(e + t) \ln (e + t)]$. From (\ref{eq:qge-two}), after multiplying by $\theta$ and integrating 

\be \label{eq:base-equation}
\frac{d}{dt} \int _{\R ^2} |\hat{\theta}|^2 \, d \xi = - 2 \int _{\R ^2} | \xi | ^{2 \a} |\hat{\theta}|^2 \, d \xi.
\ee
Then, as 

\begin{displaymath}
2 \int _{\R^2} |\xi| ^{2 \a} | \hat{\theta} | ^2 d \xi  \geq  2 \int _{B(t)} |\xi| ^{2 \a} | \hat{\theta} | ^2 d \xi + (1 + \frac{1}{\a}) [(e + t) \ln (e + t)]^{-1} \int _{B(t) ^c} | \hat{\theta} | ^2 d \xi 
\end{displaymath}
(\ref{eq:base-equation}) becomes

\bga
\frac{d}{dt} \int _{\R ^2} |\hat{\theta}|^2 \, d \xi & + & (1 + \frac{1}{\a}) [(e + t) \ln (e + t)]^{-1} \int _{\R ^2} | \hat{\theta} | ^2 d \xi \nonumber \\ & \leq & (1 + \frac{1}{\a}) [(e + t) \ln (e + t)]^{-1} \int _{B(t)} | \hat{\theta} | ^2 d \xi. \nonumber
\eda
Multiplying on both sides by $h(t) = [\ln (e+t)] ^{1 + \frac{1}{\a}}$, writing the left hand side as a derivative and integrating between $0$ and $t$ 

\bga \label{eq:ineq-prior-to-computation}
[\ln (e+t)] ^{1 + \frac{1}{\a}} \int _{\R ^2} |\hat{\theta}|^2 \, d \xi & \leq & \| \theta_0 \| _{\eledos} ^2 \nonumber \\ & + & \int _0 ^t (1 + \frac{1}{\a}) [(e + s)]^{-1}  \ln (e + s) ^{\frac{1}{\a}} \left( \int _{B(s)}  |\hat{\theta}|^2 \, d \xi \right) \, ds. 
\eda
Hence, we need to estimate $\int _{B(s)}  |\hat{\theta}|^2 \, d \xi$. From the solution to the Fourier transform of (\ref{eq:qge-two})

\begin{displaymath}
\hat{\theta}(\xi,t) = \fourierzero (\xi) e ^{- |\xi|^{2\a}t} + \int _0 ^t e ^{- |\xi|^{2\a}(t-s)} \widehat{\nonlinear} (s) ds.
\end{displaymath}
we obtain

\begin{displaymath}
|\hat{\theta}(\xi,t)| \leq |\fourierzero(\xi)| + \int _0 ^t |\widehat{\nonlinear}| \, ds 
\end{displaymath}
which, by Lemma \ref{lm:lemma-2} leads to

\bga
|\hat{\theta}(\xi,t)| ^2 & \leq & 2 \left( |\fourierzero(\xi)|^2 + \left( \int _0 ^t |\xi| \| \theta({\tau}) \| _{\eledos} ^2 \, d \tau \right) ^2 \right) \nonumber \\ & \leq & 2 \left( |\fourierzero(\xi)|^2 +  t |\xi| ^2 \int _0 ^t \| \theta({\tau}) \| _{\eledos} ^4 \, d \tau \right) \nonumber.
\eda
Then, passing to polar coordinates

\bga \label{eq:expression-for-bs}
\int _{B(s)} | \hat{\theta} | ^2 d \xi & \leq & 2 \left( \int _{B(s)} |\fourierzero|^2 \, d \xi +  \int _{B(s)} s |\xi| ^2 \left( \int _0 ^s \| \theta({\tau}) \| _{\eledos} ^4 \, d \tau \right) \, d \xi \right) \nonumber \\ & \leq & 2 \left( C [(e + s) \ln (e + s)] ^{-\frac{1}{\a}(\frac{2}{p} - 1)} + \int _0 ^{\frac{\pi}{2}} \int _0 ^{g(s)^{\frac{-2}{\a}}} r^3 \, s \left( \int _0 ^s \| \theta({\tau}) \| _{\eledos} ^4 \, d \tau \right) \, dr \, d \varphi \right) \nonumber \\ & \leq & 2 \left( C [(e + s) \ln (e + s)] ^{-\frac{1}{\a}(\frac{2}{p} - 1)} + C s^2 g^{\frac{2}{\a}}(s) \right) \nonumber \\ & = &  C [(e + s) \ln (e + s)] ^{-\frac{1}{\a}(\frac{2}{p} - 1)} + C s^2 [(e + s) \ln (e + s)]^{-\frac{2}{\a}}
\eda
where we have used Lemma \ref{lm:lemma-1} with $g(s) = C [(e + s) \ln (e + s)]$ and the Maximum Principle for the $L^2$ norm of $\theta$. Substituting (\ref{eq:expression-for-bs}) in (\ref{eq:ineq-prior-to-computation}) we see that the integral in the right hand side of  (\ref{eq:ineq-prior-to-computation}) is finite, so 

\begin{displaymath}
\int _{\R^2} | \hat{\theta} | ^2 d \xi \leq C \, [\ln (e + t)]^{-(1 + \frac{1}{\a})}.
\end{displaymath}
The formal part of the proof is now complete. To extend the estimate to weak solutions, we repeat the argument, applying it to the solutions of the approximate equations

\begin{displaymath}
\frac{\partial \theta_n}{\partial t} + u_n \nabla \theta _n + (-\Delta)^{\a} \theta _n = 0
\end{displaymath}
where $u_n =  \Psi _{\delta_n} (\theta _n)$ is defined by

\begin{displaymath}
\Psi _{\delta_n} (\theta _n) = \int _0 ^t \phi (\tau) {\mathcal R}^{\bot} \theta _n (t - \delta_n \tau) \, d \tau.
\end{displaymath} 
Here the operator ${\mathcal R}^{\bot}$ is defined on scalar functions as

\begin{displaymath}
{\mathcal R}^{\bot} f = (- \partial _{x_2} \Lambda ^{-1}f, \partial _{x_1} \Lambda ^{-1}f)
\end{displaymath}
and $\phi$ is a smooth function with support in $[1,2]$ and such that $\int_0 ^{\infty} \phi(t) \, dt = 1$. For each $n$, the values of $u_n$ depend only on the values of $\theta _n$ in $[t - 2 \delta _n, t - \delta _n]$. As stated in Constantin and Wu \cite{constantin-wu}, the functions $\theta _n$ converge to a weak solution $\theta$ and strongly in $L^2$ almost everywhere. Since the estimates obtained do not depend on $n$, they are valid for the limit function $\theta$. The proof is now complete. $\square$

\subsection{\pf  of Theorem \ref{decay-initial-l2lp}} As before, we prove a formal estimate for smooth solutions, which then can  be extended to weak solutions by the method of retarded mollifiers. 
For the formal estimate, we proceed as in Lemma \ref{lm:preliminary-estimate}, employing the Fourier splitting method with 

\begin{displaymath}
B(t) = \{ \xi: |\xi| \leq g(t) ^{-\frac{1}{2\a}} \}
\end{displaymath}
for $g(t) = 2 \a (t+1)$. Thus 

\begin{displaymath}
\frac{d}{dt} \int _{\R ^2} |\hat{\theta}|^2 \, d \xi + \frac{1}{\a (t+1)} \int _{\R ^2} |\hat{\theta}|^2 \, d \xi \leq \frac{1}{\a (t+1)} \int _{B(t)} |\hat{\theta}|^2 \, d \xi
\end{displaymath}
which after using $h(t) = (t+1)^{\frac{1}{\a}}$ as an integrating factor leads to

\be \label{eq:important-ineq-real-decay}
(t+1)^{\frac{1}{\a}} \int _{\R ^2} |\hat{\theta}|^2 \, d \xi \leq \| \fourierzero \| _{\eledos} ^2 + \int _0 ^t \frac{1}{\a} (s + 1) ^{\frac{1}{\a} -1} \left( \int _{B(s)} |\hat{\theta}|^2 \, d \xi \right) \, ds.
\ee
Working as in (\ref{eq:expression-for-bs}) in Lemma \ref{lm:preliminary-estimate}, using Lemma \ref{lm:lemma-1} with $g(s) = 2 (s+1)$ and the preliminary estimate from Lemma \ref{lm:preliminary-estimate} for making 

\begin{displaymath}
\Vert \theta (\tau) \Vert _{\eledos} ^4 \leq  \Vert \theta (\tau) \Vert _{\eledos} ^2 [\ln (e + \tau)]^{-(1 + \frac{1}{\a})}
\end{displaymath}
we obtain

\bga \label{eq:long-form}
(t+1)^{\frac{1}{\a}} \int _{\R ^2} |\hat{\theta}|^2 \, d \xi & \leq & \| \fourierzero \| _{\eledos} ^2 + C \int _0 ^t (1 + s) ^{-\frac{1}{\a}(\frac{2}{p} - 1) + \frac{1}{\a} - 1} ds \nonumber \\ & + & \int _0 ^t \int _0 ^s \| \theta({\tau}) \| _{\eledos} ^2 [\ln (e + \tau)]^{-(1 + \frac{1}{\a})} s g^{\frac{-2}{\a}}(s) (1 + s) ^{\frac{1}{\a} - 1} \, d \tau \, ds \nonumber \\ & \leq & \| \fourierzero \| _{\eledos} ^2 + C (1 + s) ^{-\frac{1}{\a}(\frac{2}{p} - 1) + \frac{1}{\a}} \nonumber \\ & + & \int _0 ^t \int _0 ^s \| \theta({\tau}) \| _{\eledos} ^2 [\ln (e + \tau)]^{-(1 + \frac{1}{\a})} s  (1 + s) ^{-(\frac{1}{\a} + 1)} \, d \tau \, ds.
\eda
Now

\bga
I(t) & = & \int _0 ^t \int _0 ^s \| \theta(\tau) \| _{\eledos} ^2 [\ln (e + \tau)]^{-(1 + \frac{1}{\a})} s  (1 + s) ^{-(\frac{1}{\a} + 1)} \, d \tau \, ds \nonumber \\ & \leq & C \int _0 ^ t (1 + s) ^{-\frac{1}{a}} \, ds \, \int _0 ^t (1 + \tau) ^{\frac{1}{\a}} \| \theta(\tau) \| _{\eledos} ^2 \frac{[\ln (e + \tau)]^{-(1 + \frac{1}{\a})}}{(1 + \tau) ^{\frac{1}{\a}}} \, d \tau \nonumber \\ & \leq &  C \int _0 ^t (1 + \tau) ^{\frac{1}{\a}} \| \theta(\tau) \| _{\eledos} ^2 \frac{[\ln (e + \tau)]^{-(1 + \frac{1}{\a})}}{(1 + \tau) ^{\frac{1}{\a}}} \, d \tau \nonumber
\eda
so taking 

\begin{displaymath}
f(t) = (t+1)^{\frac{1}{\a}} \| \theta(t) \| _{\eledos} ^2, \quad a(t) = (1 + t) ^{-\frac{1}{\a}(\frac{2}{p} - 1) + \frac{1}{\a}}, \quad b(t) = \frac{[\ln (e + \tau)]^{-(1 + \frac{1}{\a})}}{(1 + \tau) ^{\frac{1}{\a}}}
\end{displaymath}
equation (\ref{eq:long-form}) becomes

\begin{displaymath}
f(t) \leq C + a(t) + \int _0 ^t f(\tau) b(\tau) d \tau.
\end{displaymath}
By Gronwall's inequality

\be \label{eq:gronwall}
f(t)  \leq f(0) \exp \left( \int _0 ^t b(\tau) \, d \tau \right) + \int _0 ^ t a'(\tau) \exp \left( \int _\tau ^t b(s) ds \right) d \tau.
\ee
Notice that as $\frac{1}{2} < \a \leq 1$

\begin{displaymath}
\int _0 ^t b(\tau) \, d \tau = \int _0 ^t \frac{[\ln (e + \tau)]^{-(1 + \frac{1}{\a})}}{(1 + \tau) ^{\frac{1}{\a}}} < \infty. 
\end{displaymath}
Then (\ref{eq:gronwall}) becomes 

\begin{displaymath}
(t+1)^{\frac{1}{\a}} \| \hattheta (t) \| _{\eledos} ^2 \leq C \, \| \fourierzero \| _{\eledos} ^2 +  (1 + t) ^{-\frac{1}{\a}(\frac{2}{p} - 1) + \frac{1}{\a}} ds
\end{displaymath}
hence

\begin{displaymath}
\| \theta(t) \| _{\eledos} ^2 \leq (t+1)^{- \frac{1}{\a}} + C(1 + t) ^{\frac{1}{\a}(\frac{2}{p} - 1)} \leq C(1 + t) ^{\frac{1}{\a}(\frac{2}{p} - 1)}
\end{displaymath}
which proves the formal estimate. The retarded mollifiers method allows us to extend it to weak solutions. $\square$

\section{$L^q$ decay, for $q \geq \frac{2}{2\a - 1}$}
\label{proof-of-remaining-theorems}

\subsection{Proof of Theorem \ref{kato-like}} We now describe the main ideas behind the proof of Theorem \ref{kato-like}. For clarity, we let $m = \exponent$. We first prove preliminary estimates of the form

\bga
\label{eq:first-preliminary-estimate} \Vert t^{\frac{1}{\a}(\frac{1 - \delta}{m})} \theta (t) \Vert _{\frac{m}{\delta}} \leq C, \quad t > 0 \\
\label{eq:second-preliminary-estimate} \Vert t ^{\frac{1}{2 \a}} \nabla \theta (t) \Vert _{L^m} \leq C, \quad t > 0
\eda
for fixed $0 < \delta < 1$. To do so, following Katos's \cite{kato} ideas, we construct a solution in $\eledosdelta$ to the integral equation (\ref{eq:int-form}) by succesive approximations

\begin{displaymath}
\theta _1 (t) = \kernel (t) * \theta _0
\end{displaymath}

\begin{displaymath}
\theta_{n+1} (t) =  \kernel (t) * \theta _0 - \int _0 ^t \kernel (t-s) * \nonlinearapprox (s) \, ds, \quad n \geq 1.
\end{displaymath}
These approximations are such that

\be \label{eq:iteration-step}
\Vert t^{\frac{1}{\a}(\frac{1 - \delta}{m})} \theta_{n + 1} (t) \Vert _{L^{\frac{m}{\delta}}} \leq K_{n+1}, \quad \Vert  t ^{\frac{1}{2 \a}}  \nabla \theta _{n+1} (t) \Vert _{L^m} \leq K' _{n+1}
\ee
are bounded by expressions that depend on $K_1, K_n$ and $K' _n$ only. If $\theta _0$ has small $L^m$ norm then these recursive relations are uniformly bounded, this is 

\begin{displaymath}
K_n \leq K, \, K' _n \leq K, \quad n \geq 1
\end{displaymath}
for some $K > 0$. A standard argument allows us to show that there is a uniformly converging subsequence $\theta_n$ whose limit is a solution to (\ref{eq:int-form}) in $\eledosdelta$ that obeys (\ref{eq:first-preliminary-estimate}) and (\ref{eq:second-preliminary-estimate}). These preliminary estimates are used to bootstrap a similar argument which proves the results stated in the Theorem.

\pf  We begin by proving (\ref{eq:first-preliminary-estimate}) and (\ref{eq:second-preliminary-estimate}). Let $\delta$ be fixed, $0 < \delta < 1$. We note first that by (\ref{eq:norm-ineq-one-lemma}) in Lemma \ref{prop-wu}

\begin{displaymath}
\Vert \theta _1 (t) \Vert _{\eledosdelta} \leq C  t^{- \frac{1}{\a}(\frac{1 - \delta}{m})} \Vert \theta _0 \Vert _{L^m}. 
\end{displaymath}
and by (\ref{eq:norm-ineq-two-lemma}) in Lemma \ref{prop-wu}

\begin{displaymath}
\Vert \nabla \theta _1 (t) \Vert _{L^m} \leq C t ^{- \frac{1}{2 \a}} \Vert \theta _0 \Vert _{L^m}.
\end{displaymath}
Let $K_1 = K' _1 = C \Vert \theta_0 \Vert _{L^m}$. Now assume 

\bga
\Vert t^{\frac{1}{\a}(\frac{1 - \delta}{m})} \theta_{n} (t) \Vert _{\frac{m}{\delta}} \leq K_n \nonumber \\
\Vert t ^{\frac{1}{2 \a}} \nabla \theta _{n} (t) \Vert _{L^m} \leq K' _n \nonumber
\eda
for $t > 0$. Then

\bga \label{eqn:bound-norm-thetanplusone}
\Vert \theta _{n+1} (t) \Vert _{\eledosdelta} & \leq & \| \theta_1 (t) \| _{\eledosdelta} + \int _0 ^t \Vert \kernel (t-s) * \nonlinearapprox (s) \Vert _{\eledosdelta} \, ds \nonumber \\ & \leq & K_1 \, t^{- \frac{1}{\a}(\frac{1 - \delta}{m})} + \int _0 ^t \Vert \kernel (t-s) * \nonlinearapprox (s) \Vert _{\eledosdelta} \, ds \nonumber \\ & \leq & K_1 \, t^{- \frac{1}{\a}(\frac{1 - \delta}{m})} + C \int _0 ^t (t-s) ^{-\frac{1}{\a m}} \Vert \theta_n (s) \Vert _{\eledosdelta} \Vert \nabla \theta _n (s) \Vert _{L^m} \, ds \nonumber \\ & \leq & K_1 \, t^{- \frac{1}{\a}(\frac{1 - \delta}{m})} + C \, K_n \, K' _n \int _0 ^t (t-s) ^{-\frac{1}{\a m }} s^{- \frac{1 - \delta}{\a m} - \frac{1}{2 \a}} \, ds \nonumber \\ & \leq &  K_1  \, t^{- \frac{1}{\a}(\frac{1 - \delta}{m})} +  C \, K_n \, K' _n t^{- \frac{1}{\a}(\frac{1 - \delta}{m})}
\eda
where we used boundedness of the Riesz transform and (\ref{eq:first-ineq-int-norm}) in Lemma \ref{main-lemma-ineq} with $\eta = \mu = \frac{2 \delta}{m}$ and $\nu = \frac{2}{m}$. By an analogous method

\bga \label{eqn:bound-norm-gradthetanplusone}
\Vert \nabla \theta _{n+1} (t) \Vert _{L^m} & \leq & \| \nabla \theta_1 (t) \| _{L^m} + \int _0 ^t \Vert \nabla \kernel (t-s) * \nonlinearapprox (s) \Vert _{L^m} \, ds \nonumber \\ & \leq & K_1 \, t ^{- \frac{1}{2 \a}} + \int _0 ^t \Vert \nabla \kernel (t-s) * \nonlinearapprox (s) \Vert _{L^m} \, ds \nonumber \\ & \leq & K_1 \, t ^{- \frac{1}{2 \a}}  + C \int _0 ^t (t-s) ^{-(\frac{1}{2 \a}+ \frac{\delta}{m \a})} \Vert \theta_n (s) \Vert _{L^{\frac{m}{delta}}} \Vert \nabla \theta _n (s) \Vert _{L^m} \, ds \nonumber \\ & \leq & K_1 \, t ^{- \frac{1}{2 \a}} + C \, K_n \, K' _n \int _0 ^t (t-s) ^{-(\frac{1}{2 \a} + \frac{\delta}{\a m})} s^{- \frac{1 - \delta}{\a m} - \frac{1}{2 \a}} \, ds \nonumber \\ & \leq &  K_1 \, t ^{- \frac{1}{2 \a}}  +  C \, K_n \, K' _n t ^{- \frac{1}{2 \a}}
\eda
where we used  boundedness of the Riesz transform and (\ref{eq:second-ineq-int-norm}) in Lemma \ref{main-lemma-ineq} with $\eta = \nu = \frac{2}{m}, \mu = \frac{2 \delta}{m}$. We have then that the norms described in (\ref{eq:iteration-step}) are respectively bounded by 

\bga
K_{n+1} \leq K_1 + C K_n K'_n \nonumber \\
K' _{n+1} \leq K_1 + C K_n K'_n. \nonumber
\eda
If $K_1 < \frac{1}{4c}$, an  induction argument allows us to prove that

\begin{displaymath}
K_n \leq K, \quad K' _n \leq K
\end{displaymath}
for $n \geq 1$, where $K = \frac{1}{2c}$. Note that $K_0 < \frac{1}{4c}$ implies $\Vert \theta_0 \Vert _{L^m} < \frac{1}{4c^2}$, thus the $L^m$ norm of the initial data has to be small. Then

\bga
\Vert t^{\frac{1}{\a}(\frac{1 - \delta}{m})} \theta_{n} (t) \Vert _{L^{\frac{m}{\delta}}} \leq K \nonumber \\
\Vert  t ^{\frac{1}{2 \a}}  \nabla \theta _{n} (t) \Vert _{L^m} \leq K. \nonumber
\eda
for $n \geq 1$. By a standard argument (see Kato \cite{kato} and Kato and Fujita \cite{kato-fujita} for full details) we can extract a subsequence that converges uniformly in $(0, +\infty)$ to a solution $\theta$. Then

\bga
t^{\frac{1}{\a}(\frac{1 - \delta}{m})} \theta \in BC((0, + \infty), \eledosdelta) \nonumber \\
t ^{\frac{1}{2 \a}}  \nabla \theta \in BC((0, +\infty), L^m) \nonumber.
\eda

We now use these preliminary estimates to prove the Theorem. As before, we construct a solution by sucessive approximations. Let $m \leq q < \infty$. By (\ref{eq:norm-ineq-one-lemma}) in Lemma \ref{prop-wu}

\begin{displaymath}
\Vert \theta _1  (t) \Vert _{\eleq} \leq C t^{-\frac{1}{\a}(\frac{1}{m} - \frac{1}{q})} \Vert \theta _0 \Vert _{L^m}
\end{displaymath}
and by (\ref{eq:norm-ineq-two-lemma}) in Lemma \ref{prop-wu}

\begin{displaymath}
\Vert \nabla \theta _1 (t) \Vert _{\eleq} \leq C t^{-(\frac{1}{2 \a} + \frac{1}{\a}(\frac{1}{m} - \frac{1}{q}))} \Vert \theta_0 \Vert _{L^m} .
\end{displaymath}
Notice that this estimate holds for $q \geq m$. Again, set $K_1 = K' _1 = C \Vert \theta_0 \Vert _{L^m}$. We want to show inductively that the $L^q$ norms of $t^ {\frac{1}{\a}(\frac{1}{m} - \frac{1}{q})} \th _n (t)$ and $t^{\frac{1}{2 \a}+ \frac{1}{\a}(\frac{1}{m} - \frac{1}{q})} \nabla_{n + 1} \th (t)$ are uniformly bounded. Then
 
\bga
\Vert \theta_{n+1} (t)  \Vert _{\eleq} & \leq & K_1 t^{-\frac{1}{\a}(\frac{1}{m} - \frac{1}{q})} +  \int _0 ^t \Vert \kernel (t-s) * \nonlinearapprox (s) \Vert _{\eleq} \, ds \nonumber \\ & \leq & K_1 t^{-\frac{1}{\a}(\frac{1}{m} - \frac{1}{q})} + C \int _0 ^t (t-s) ^{-\frac{1}{\a}(\frac{1 + \delta}{m} - \frac{1}{q})} \Vert \theta_n (s) \Vert _{\eledosdelta} \Vert \nabla \theta_n (s) \Vert _{L^m} \, ds \nonumber \\ & \leq & K_1 t^{-\frac{1}{\a}(\frac{1}{m} - \frac{1}{q})} + C K_n K '_n \int _0 ^t (t-s) ^{-\frac{1}{\a}(\frac{1 + \delta}{m} - \frac{1}{q})}  s^{- \frac{2 - \delta}{m \a} - \frac{1}{2 \a}} \, ds \nonumber \\ & \leq & K_1 \, t^{-\frac{1}{\a}(\frac{1}{m} - \frac{1}{q})} + C K_n K'_n t^{-\frac{1}{\a}(\frac{1}{m} - \frac{1}{q})}  
\eda
where we have used (\ref{eq:first-ineq-int-norm}) in Lemma \ref{main-lemma-ineq} with $\eta = \frac{2}{q}, \mu = \frac{2 \delta}{m}$ and $\nu = \frac{2}{m}$ and we have used the preliminary estimates obtained for $\| \theta _n (t) \| _{\eledosdelta}$ and $\| \nabla \theta _n (t) \| _{L^m}$. Proceeding analogously for the gradient we obtain

\bga
\Vert \nabla \theta_{n+1} (t)  \Vert _{\eleq} & \leq & \| \nabla \theta _1 (t) \| _{\eleq} + \int _0 ^t \Vert \nabla \kernel (t-s) * \nonlinearapprox (s) \Vert _{\eleq} \, ds \nonumber \\ & \leq & K_1 t^{-(\frac{1}{2 \a}+ \frac{1}{\a}(\frac{1}{m} - \frac{1}{q}))} +  \int _0 ^t \Vert \nabla \kernel (t-s) * \nonlinearapprox (s) \Vert _{\eleq} \, ds \nonumber \\ & \leq & K_1 t^{-(\frac{1}{2 \a}+ \frac{1}{\a}(\frac{1}{m} - \frac{1}{q}))} + C \int _0 ^t (t-s) ^{-(\frac{1}{2 \a} + \frac{1}{\a}(\frac{1 + \delta}{m} - \frac{1}{q}))} s^{- \frac{2 - \delta}{m \a} - \frac{1}{2 \a}} \, ds \nonumber \\ & \leq & K_1 t^{-(\frac{1}{2 \a}+ \frac{1}{\a}(\frac{1}{m} - \frac{1}{q}))} + K_n K' _n t^{-(\frac{1}{2 \a}+ \frac{1}{\a}(\frac{1}{m} - \frac{1}{q}))}
\eda
where we used (\ref{eq:second-ineq-int-norm}) in Lemma \ref{main-lemma-ineq} with $\eta = \frac{2}{q}, \mu = \frac{2 \delta}{m}$ and $\nu = \frac{2}{m}$. As before, setting 

\bga
K_{n+1} =  \Vert t^{\frac{1}{\a}(\frac{1}{m} - \frac{1}{q})} \th _{n+1} (t) \Vert _{L^q} \nonumber \\
K' _{n+1} = \Vert t^{\frac{1}{2 \a}+ \frac{1}{\a}(\frac{1}{m} - \frac{1}{q})} \nabla \th_{n+1} (t) \Vert _{L^q}
\eda
we obtain  

\bga
K_{n+1} \leq K_1 + C K_n K'_n \nonumber \\
K' _{n+1} \leq K_1 + C K_n K'_n. \nonumber
\eda
The same arguments that were used for the preliminary estimates apply here, leading to 

\begin{displaymath}
t^ {\frac{1}{\a}(\frac{1}{m} - \frac{1}{q})} \th (t) \in BC((0, \infty), \eleq)
\end{displaymath}
and

\begin{displaymath}
t^{\frac{1}{2 \a}+ \frac{1}{\a}(\frac{1}{m} - \frac{1}{q})} \nabla \th (t) \in BC ((0, \infty), \eleq)
\end{displaymath}
for $ \frac{2}{2 \a -1} \leq q < \infty$, which is the desired result. $\square$ 
\subsection{Proof of Theorem \ref{theo-decay-2}} Let $m = \frac{2}{2\a - 1}$. By (\ref{eq:max-prin-ju}), the $L^m$ norm of $\theta$ tends to zero, so for times larger than some $T = T(\theta_0)$, $\Vert \theta (t) \Vert _{L^m} \leq \kappa$, for $\kappa$ as in Theorem \ref{kato-like}. Let $m \leq q < r$. Interpolation yields 

\begin{displaymath}
\| \theta (t) \| _{\eleq} \leq \| \theta (t) \| _{L^m} ^a \, \| \theta (t) \| _{L^r} ^{1 - a}
\end{displaymath}
for $a = \frac{m}{q} \frac{r - q}{r - m}$ and $1- a = \frac{r}{q} \frac{q - m}{r - m}$. Then

\begin{displaymath}
\| \theta (t) \| _{\eleq} \leq  C t^{(1 - \frac{1}{\a})\frac{m}{q} \left( \frac{r - q}{r - m} \right) - \frac{1}{\a}(\frac{1}{m} - \frac{1}{q})}.
\end{displaymath}
This holds for any $r$ such that $q \leq r < \infty$. The optimal decay rate is given by the minimum of the exponent

\begin{displaymath}
f(r) = C_1 \, \frac{r - q}{r - m} - C_2
\end{displaymath} 
where  $C_1 = (1 - \frac{1}{\a})\frac{m}{q}$ and $C_2 = \frac{1}{\a}(\frac{1}{m} - \frac{1}{q})$. As $C_1 < 0$, this is a non-increasing function, so the optimal decay rate is 

\begin{displaymath}
\lim _{r \to \infty} f(r) = C_1 - C_2 = \left( 1 - \frac{1}{\a} \right) \frac{m}{q} - \frac{1}{\a} \left( \frac{1}{m} - \frac{1}{q} \right).
\end{displaymath}
Then

\begin{displaymath}
\| \theta (t) \| _{L^q} \leq C \, t^{\frac{1}{q}\frac{4a - 3}{\a(2\a - 1)} - 1 + \frac{1}{2 \a}}, \quad t \geq T
\end{displaymath}
 which is the desired result. $\square$

\end{document}